\documentclass[10pt]{amsart}
\usepackage{amsmath, amssymb, amsthm, a4wide}
\usepackage{geometry}
\usepackage{fancyhdr}
\usepackage{array}
\usepackage{enumerate}
\newtheorem{theorem}{Theorem}[section]
\newtheorem{lemma}[theorem]{Lemma}
\newtheorem{definition}[theorem]{Definition}

\newtheorem{corollary}[theorem]{Corollary}

\newtheorem{remark}[theorem]{Remark}

\numberwithin{equation}{section}

\newcounter{newlist}
\newenvironment{mylist}[1][]{
\begin{list}{\textnormal{(\arabic {newlist})}}
    {
    \usecounter{newlist}
     \setlength{\labelsep}{0.5em}
     \setlength{\leftmargin}{0.5em}
     \setlength{\rightmargin}{0cm}
     \setlength{\parsep}{0pt}
     \setlength{\itemsep}{3pt}
     \setlength{\itemindent}{0em}
     \setlength{\listparindent}{2em}
    }}
{\end{list}}
\newcounter{nnnewlist}

\newcounter{nelist}
\newenvironment{mlist}[1][]{
\begin{list}{\textnormal{(A\arabic {nelist})}}
    {
    \usecounter{nelist}
     \setlength{\labelsep}{0.5em}
     \setlength{\leftmargin}{2.5em}
     \setlength{\rightmargin}{0cm}
     \setlength{\parsep}{0pt}
     \setlength{\itemsep}{3pt}
     \setlength{\itemindent}{0em}
     \setlength{\listparindent}{2em}
    }}
{\end{list}}
\newcounter{mnelist}
\newenvironment{mnlist}[1][]{
\begin{list}{\textnormal{(A\arabic {mnelist}')}}
    {
    \usecounter{mnelist}
     \setlength{\labelsep}{0.5em}
     \setlength{\leftmargin}{2.5em}
     \setlength{\rightmargin}{0cm}
     \setlength{\parsep}{0pt}
     \setlength{\itemsep}{3pt}
     \setlength{\itemindent}{0em}
     \setlength{\listparindent}{2em}
    }}
{\end{list}}

\begin{document}

\title[RGSDEs with nonlinear resistance]{Reflected Stochastic Differential Equations\\ Driven by $G$-Brownian Motion\\ with Nonlinear Resistance}
\author{Peng Luo}
\thanks{School of Mathematics and Qilu Securities Institute for Financial Studies, Shandong University and Department of Mathematics and Statistics, University of Konstanz; pengluo1989@gmail.com}
\thanks{The author's research was
partially supported by China Scholarship Council, NSF (No. 10921101) and by the 111Project (No. B12023).}

%\subjclass[2000]{60H10}

\date{}


\begin{abstract}
In this paper, we study the uniqueness and existence of solutions of RGSDEs with nonlinear resistance under an integral-Lipschitz condition of coefficients. Moreover we obtain the comparison theorem for RGSDEs with nonlinear resistance.
\end{abstract}

\maketitle

\section{Introduction}
\noindent In classical framework,  the diffusion processes with reflecting boundaries were introduced by Skorokhod \cite{S1,S2} in 1960s. After that, many works related to reflected solutions to SDEs and BSDEs have been done. El Karoui \cite{EL}, El Karoui and Chaleyat-Maurel \cite{EC} and Yamada \cite{Y} studied scalar valued reflected SDEs on a half-line. For multidimensional case, Stroock and Varadhan \cite{SV} obtained the existence of weak solutions to reflected SDEs on a smooth domain which was extended to a convex domain by Tanaka \cite{T} and  a non convex domain by Lions and Sznitman \cite{LS}. On the other hand, the solvability of reflected BSDEs was first obtained by El Karoui et al.\cite{EKPPQ}. Then many corresponding results for reflected BSDEs have been established by Gegout-Petit and Pardoux \cite{GP}, Ramasubramanian \cite{R} and Hu and Tang \cite{HT}, etc.. In particular, Qian and Xu \cite{QX} obtained the existence and uniqueness of solutions of reflected BSDEs with nonlinear resistance under a Lipschitz condition.\\[6pt]
Motivated by uncertainty problems, risk measures and the superhedging in finance, Peng systemically established a time-consistent fully nonlinear
expectation theory (see \cite{Peng2004,Peng2005,P09}). As a typical and important case, Peng  introduced the
$G$-expectation theory (see \cite{Peng4, P10}  and the references therein) in 2006.
In the  $G$-expectation framework ($G$-framework for short), the
notion of  $G$-Brownian motion and the corresponding stochastic
calculus of It\^{o}'s type were established. On that basis, Gao \cite{G} and Peng \cite {Peng4}  studied the existence and uniqueness of the solution of $G$-SDE under a standard Lipschitz condition. For a recent account and development of this theory we refer the reader to \cite{L1, LQ1, LQ, LP, LY,BL, LW, LJ, P2, XZ}. Recently, Lin \cite{L} obtained the existence and uniqueness of the solution of $G$-SDE with reflecting boundary. Compared to Lin \cite{L}, we consider an integral-Lipschitz condition of coefficients and also the increasing process $K$ contributes to the coefficients, namely the following scalar valued RGSDE with nonlinear resistance:
\begin{equation}\label{ABC}
\left\{\begin{aligned}
&X_t=x+\int^t_0 f_s(X_s,K_s)ds+\int^t_0 h_s(X_s,K_s)d\langle B \rangle_s +\int^t_0 g_s(X_s,K_s) dB_s +K_t, \ q.s.,\ 0\leq t\leq T; \\
&X_t\geq S_t; \int^T_0 (X_t-S_t)dK_t=0,
\end{aligned}\right.
\end{equation}
where $\langle B\rangle$ is the quadratic variation process of $G$-Brownian motion $B$, and $K$ is an increasing process which pushes the solution $X$ upwards to be remaining above the obstacle $S$ in a minimal way. The aim of this paper is to study the existence and uniqueness of solutions to the above RGSDEs with nonlinear resistance in the sense of ``quasi-surely'' defined by Denis et al. \cite{DHP}. The main idea is to estimate in some sense simultaneously the solution $X$ and the increasing process $K$ from which the uniqueness result follows and  a solution in $M^p_G([0, T])$ to (\ref{ABC}) can be constructed by fixed-point iteration. To establish the comparison theorem, we use the extented $G$-It\^{o} formula in Lin \cite{L}.\\[6pt]
This paper is organized as follows: Section 2 introduces some notations and results in $G$-framework which is necessary for what follows, while section 3 is our main results.
\section{Preliminaries}
The main purpose of this section is to recall some preliminary
results in $G$-framework which are needed in the sequel.  More
details can be found in  Denis et al \cite{DHP}, Li and Peng \cite{LP}, Lin \cite{L,LY} and Peng \cite{Peng4}.

Denote by
$\Omega=C_{0}^{d}(\mathbb{R}^{+})$ the space of all
$\mathbb{R}^{d}$-valued continuous paths
$(\omega_{t})_{t\in\mathbb{R}^{+}}$, with $\omega_{0}=0$, equipped
with the distance
$$\rho(\omega^{1},\omega^{2}):=\sum_{i=1}^{\infty}2^{-i}[\max\limits_{t\in[0,i]}|\omega^{1}_t-\omega^{2}_t|\wedge 1].$$
$\mathcal{B}({\Omega})$ is the Borel $\sigma$-algebra of $\Omega$.

For each $t\in[0,\infty)$, we introduce the following spaces.
\\
$\bullet$ $\Omega_{t}:=\{\omega({\cdot\wedge t}):\omega\in\Omega\}$, $\mathcal{F}_{t}:=\mathcal {B}(\Omega_{t})$,
\\
$\bullet$ $L^{0}{(\Omega)}:$ the space of all $\mathcal
{B}(\Omega)$-measurable real functions,
\\
$\bullet$ $L^{0}{(\Omega_{t})}:$ the space of all $\mathcal{F}_{t}$-measurable real functions,
\\
$\bullet$ $B_{b}{(\Omega)}:$ all bounded elements in $L^{0}{(\Omega)}$, $B_{b}{(\Omega_{t})}:=B_{b}{(\Omega)}\cap L^{0}{(\Omega_{t})}$,
\\
$\bullet$
$C_{b}{(\Omega)}:$ all  continuous elements in $B_{b}{(\Omega)}$, $C_{b}{(\Omega_{t})}:=C_{b}{(\Omega)}\cap L^{0}{(\Omega_{t})}.$

In Peng \cite{Peng4}, a $G$-Brownian motion is constructed on a sublinear expectation
space $(\Omega,{L}_{G}^1,\hat{\mathbb{E}},(\hat{\mathbb{E}}_t)_{t\geq 0})$, where   ${L}_{G}^p(\Omega)$
 is a Banach space
under the natural norm $\|X\|_p=\hat{\mathbb{E}}[|X|^p]^{1/p}$.   In this space the corresponding
canonical process $B_t(\omega) = \omega_t$
 is a $G$-Brownian motion. Denote $L^p_b(\Omega)$ the completion of $B_b(\Omega)$.
Denis et al.\cite{DHP}  proved that
$L^{p}_b{(\Omega)}\supset{L}_{G}^{p}(\Omega)\supset
C_{b}{(\Omega)}$, and there exists a weakly compact family $\mathcal
{P}$ of probability measures defined on $(\Omega, \mathcal
{B}(\Omega))$ such that
$$\hat{\mathbb{E}}[X]=\sup\limits_{P\in\mathcal{P}}E_{P}[X], \ \ X\in{L}_{G}^{1}{(\Omega)}.$$

\begin{remark}{\upshape Denis et al. \cite{DHP} gave a concrete set $\mathcal{P}_M$ that represents $\hat{\mathbb{E}}$.
Consider a 1-dimensional Brownian motion $ {B_t}$ on $(\Omega,\mathcal{F},P)$, then
\[
\mathcal{P}_M := \{P_{\theta} : P_{\theta}= P\circ X^{-1},\ X_t = \int^t_0 \theta_sdB_s,\  \theta\in L^2_{\mathcal{F}}([0, T ]; [\underline{\sigma}^2, \overline{\sigma}^2])\}\]
is a set that represents $\hat{\mathbb{E}}$, where $L^2_{\mathcal{F}}([0, T ]; [\underline{\sigma}^2, \overline{\sigma}^2])$ is the collection of all $\mathcal{F}$-adapted
measurable processes with $\underline{\sigma}^2 \leq |\theta(s)|^2 \leq\overline{\sigma}^2$.
}
\end{remark}

Now we introduce the natural Choquet capacity
$$c(A):=\sup\limits_{P\in\mathcal{P}}P(A), \ \ A\in\mathcal
{B}(\Omega).$$
\begin{definition}A set $A\subset\mathcal{B}(\Omega)$ is polar if $c(A)=0$.  A
property holds $``quasi$-$surely''$ (q.s.) if it holds outside a
polar set.
\end{definition}
%\begin{definition}A real function $X$ on $\Omega$ is said to be quasi-continuous if for each $\varepsilon>0$,
%there exists an open set $O$ with $c(O)<\varepsilon$ such that
%$X|_{O^{c}}$ is continuous.
%\end{definition}
%
%\begin{definition}  We say that  $X:\Omega\mapsto\mathbb{R}$ has a quasi-continuous
%version if there exists a quasi-continuous function $Y:\Omega\mapsto\mathbb{R}$ such
%that $X = Y$, q.s..
%\end{definition}
%
% Then ${L}_{b}^{p}(\Omega)$ and $\mathbb{L}_{G}^{p}(\Omega)$ can be characterized as
%follows:
%$${L}_{b}^{p}(\Omega)=\{X\in {L}^{0}(\Omega)|\lim\limits_{N\rightarrow\infty}\mathbb{\hat{E}}[|X|^pI_{|X|\geq N}]=0\}$$
%and
%$$\mathbb{L}_{G}^{p}(\Omega)=\{X\in {L}^{p}_b(\Omega)|\ \  X\ \text {has a quasi-continuous version}\}.$$
Let $T\in\mathbb{R}^{+}$ be fixed.
\begin{definition}\label{caocao} For each  $p\geq 1$,
consider the following simple type of processes:
\begin{align*}
M_{G}^{0,p}([0,T])=&\{\eta:=\eta_t(\omega)=\sum_{j=0}^{N-1}\xi_{j}(\omega)I_{[t_{j},t_{j+1})}(t)\\
 &\forall\ N>0,\ 0=t_{0}<...<t_{N}=T,\  \xi_{j}\in \mathbb{L}_{G}^p(\Omega_{t_{j}}),\
 j=0,1,2,...,N-1\}.
\end{align*}
Denote by $M_{G}^{p}([0,T])$ the completion of
$M_{G}^{0,p}([0,T])$ under the norm
$$||\eta||_{M^{p}_{G}([0,T])}=|\int_{0}^{T}\hat{\mathbb{E}}[|\eta(t)|^{p}]dt|^{1/p}.$$
\end{definition}

Unlike the classical theory, the quadratic variation process of
$G$-Brownian motion $B$ is not always a deterministic process and it
can be formulated in ${L}^2_G(\Omega_{t})$ by \[\langle
B\rangle_t : =\lim\limits_{N\rightarrow\infty}\sum\limits_{i=0}^{N-1}(B_{t_{i+1}^N}-B_{t^N_i})^2= B^2_t-2 \int^t_0 B_sdB_s,\]
where $t_i^N=\frac{iT}{N}$ for each integer $N\geq 1$.

Peng \cite {Peng4} also introduced the related stochastic
calculus of It\^{o}'s type with respect to $G$-Brownian motion and the quadratic variation process (see also Li and Peng \cite{LP}), i.e., $(\int_{0}^{t}\eta_sdB_s)_{0\leq t\leq T}$ and $(\int_{0}^{t}\xi_sd\langle B\rangle_s)_{0\leq t\leq T}$ are well defined for each $\eta\in M^{2}_{G}([0,T])$ and $\xi\in M^{1}_{G}([0,T])$.

 In view of the dual formulation of
$G$-expectation as well as the properties of the quadratic variation
process $\langle B\rangle$ in $G$-framework, Gao \cite{G} obtained the following BDG type
inequalities. \begin{lemma}  For each $p \geq 1$ and $\eta \in
M^p_G([0, T])$,   \[\mathbb{\hat{E}}[ \sup\limits_ {0 \leq t\leq
T} |\int^t_0\eta_sd\langle B\rangle_s|^p]\leq
\bar{\sigma}^{2p}T^{p-1}\int^T_0\mathbb{\hat{E}} [|\eta_s|^p]ds.\]\end{lemma}
\begin{lemma}  Let $p \geq 2$ and $\eta \in M^p_G([0, T])$.
Then there exists some constant $C_p$ depending only on $p$ and $T$ such
that
\[\mathbb{\hat{E}}[ \sup\limits_ {0\leq t\leq T} |\int^t_0\eta_sd
B_s|^p]\leq C_p \mathbb{\hat{E}}[|\int^T_0
|\eta_s|^2ds|^{\frac{p}{2}}].\]\end{lemma}
 Recently Lin \cite{L} eatablished stochastic integrals with respect to an increasing process in a Riemann-Stieltjes way, i.e., $(\int_{0}^{t}X_sdK_s)_{0\leq t\leq T}$ is well defined for each $X\in M_c([0, T])$ and $K\in M_I([0, T])$, where\\
 $\bullet$ $M_c([0, T])$: the space of all processes $X$ whose paths are continuous in $t$ on $[0, T]$ outside a polar set $A$.\\
 $\bullet$ $M_I([0, T])$: the space of all q.s. increasing processes $K\in M_c([0, T])$. \\
 Moreover an extension of $G$-It\^{o}'s formula was also obtained. For more details, we refer the reader to Lin \cite{L},\cite{LY}. In particular, we recall the following argument.
\begin{lemma}\label{BXT}
For some $p>2$, consider a q.s. continuous $G$-It\^o process $Y$ defined in the following form \begin{equation}\label{gito}
Y_t= x+\int^t_0 f_sds+\int^t_0h_sd\langle
B\rangle_s+\int^t_0g_sdB_s,\ 0\leq t\leq T,
\end{equation}
where $f$, $h$ and $g$ are elements in $M^p_G([0, T])$. Then, there exists a unique pair of processes $(X, K)$ in $M^{p}_G([0, T])\times (M_I([0, T])\cap M^{p}_G([0, T]))$ such that
\begin{equation}\label{suchthat}
X_t=Y_t+K_t,\ q.s.,
\end{equation}
and (a) $X$ is positive; (b) $K_0=0$; and (c) $\int^T_0 X_t dK_t=0$, q.s..
\end{lemma}
\section{Scalar valued RGSDEs with nonlinear resistance}
\noindent In this section, we give the existence and uniqueness of the solutions to the scalar valued RGSDEs with nonlinear resistance under an intgral-Lipschitz condition. Moreover, a comparison theorem is obtained. Without loss of generality, we always assume $\overline{\sigma}^2=1$ in what follows. Before we move to our main results, we want to mention the so-called Bihari's inequality (cf. Bihari \cite{B}).
\begin{lemma}\label{le 2.2.1}
Let $\rho:(0,+\infty)\rightarrow(0,+\infty)$ be a continuous and increasing  function that vanishes at $0_+$ and satisfies $\int^1_0{\frac{dr}{\rho(r)}}=+\infty$. Let $u$ be a measurable and non-negative function defined on $(0,+\infty)$ satisfies
\[
u(t)\leq a+\int^t_0{\kappa(s)\rho(u(s))ds},~~t\in(0,+\infty)
\]
where $a\in\mathbb{R}^+$, and $\kappa:[0,T]\rightarrow\mathbb{R}^+$ is Lebesgue integrable. we have
\begin{mylist}
\item If $a=0$, then $u(t)=0,~~t\in(0,+\infty)$, $\lambda$-a.e.;
\item If $a>0$, we define
\[
\upsilon(t):=\int^t_{t_0}{\frac{1}{\rho(s)}ds}, ~~~t\in\mathbb{R}_+,
\]
where $t_0\in(0,+\infty)$, then
\[
u(t)\leq\upsilon^{-1}(\upsilon(a)+\int^t_0{\kappa(s)ds}).
\]
\end{mylist}
\end{lemma}
\noindent We consider the following scalar valued RGSDE with nonlinear resistance:
\begin{equation}\label{1}
X_t=x+\int^t_0 f_s(X_s,K_s)ds+\int^t_0 h_s(X_s,K_s)d\langle B \rangle_s +\int^t_0 g_s(X_s,K_s) dB_s +K_t,\ 0\leq t\leq T,
\end{equation}
where\\[-6pt]
\begin{mlist}
\item The initial condition $x\in \mathbb{R}$;
\item For some $p>2$, the coefficients $f$, $h$ and $g: \Omega\times[0,T]\times\mathbb{R}\times\mathbb{R}\rightarrow\mathbb{R}$ are given functions satisfying for each $x,y\in\mathbb{R}$, $f_\cdot(x,y)$, $h_\cdot(x,y)$, and $g_\cdot(x,y)\in M^p_G([0,T])$ and $|f_t(x,y)|^p+|h_t(x,y)|^p+|g_t(x,y)|^p\leq |\beta_1(t)|^p+\beta^{p}_{2}(|x|^p+|y|^p),$ where $\beta_1\in M^p_G([0,T])$ and $\beta_2\in\mathbb{R}_{+}$;
\item The coefficients $f$, $h$ and $g$ satisfying an integral-Lipschitz condition, i.e., for each $t\in[0, T]$ and $x, x', y, y'\in\mathbb{R}$, $|f_t(x,y)-f_t(x',y')|^p+|h_t(x,y)-h_t(x',y')|^p+|g_t(x,y)-g_t(x',y')|^p\leq \beta(t)\rho(|x-x'|^p+|y-y'|^p),$
 where $\beta:[0,T]\rightarrow\mathbb{R}^{+}$ is integrable, and $\rho~:(0,+\infty)\rightarrow(0,+\infty)$ is  continuous increasing and concave function that vanishes at $0_+$ and satisfies
 \[
 \int^1_0{\frac{dr}{\rho(r)}}=+\infty;
 \]
\item The obstacle is a $G$-It\^o process whose coefficients are all elements in $M^p_G([0,T])$,
and we shall always assume that $S_0\leq x$, q.s..
\end{mlist}
The solution of RGSDE with nonlinear resistance (\ref{1}) is a pair of processes $(X,K)$ which take values both in $\mathbb{R}$ and satisfy:\\[-6pt]
\begin{mlist}\addtocounter{nelist}{4}
\item $X\in M^p_G([0,T])$ and $X_t\geq S_t$, $0\leq t\leq T$, q.s.;
\item $K\in M_I([0, T])\cap M^p_G([0,T])$ and $K_0=0$, q.s.;
\item $\int^T_0(X_t-S_t)dK_t=0$, q.s..
\end{mlist}
Now we give our main results.
\begin{theorem}\label{Exi}
Let the assumptions (A1)-(A4) hold true, then the RGSDE (\ref{1}) admits a unique solution in $M^p_G([0,T])$.
\end{theorem}
\noindent In order to prove Theorem \ref{Exi}, we need some lemmas which give a prior estimate of the solution and estimate of variation in the solutions.
\begin{lemma}\label{abcde}
Assume that $(X, K)$ is a solution to (\ref{1}), then there exists a constant $C>0$ such that
\begin{equation*}
\begin{array}{r@{\ }l}
\hat{\mathbb{E}}[\sup_{0\leq s\leq T}|X_t|^p]+\hat{\mathbb{E}}[|K_T|^p]\leq C(|x|^p+\int^T_0\hat{\mathbb{E}}[|\beta_1(t)|^p]dt+\hat{\mathbb{E}}[\sup_{0\leq t\leq T}|S^+_t|^p]).
\end{array}
\end{equation*}
\end{lemma}
\noindent\textbf{Proof :}
Let $(X, K)$ be a pair of solution to (\ref{1}). Replacing $Y_t$ by
$x+\int^t_0 f_s(X_s,K_s)ds+\int^t_0 h_s(X_s,K_s)d\langle B \rangle_s +\int^t_0 g_s(X_s,K_s) dB_s-S_t$ and $X_t$ by $X_t-S_t$ in (\ref{suchthat}), we have the following representation of $K$ on $[0, T]$:
\begin{equation}\label{represen}
K_t=\sup_{0\leq s\leq t}\bigg(x+\int^s_0 f_u(X_u,K_u)du+\int^s_0 h_u(X_u,K_u)d\langle B \rangle_u +\int^s_0 g_u(X_u,K_u) dB_u-S_s\bigg)^-,\ q.s..
\end{equation}
As $X$ is the solution to (\ref{1}), we obtain
\begin{equation*}
\begin{array}{r@{\ }l}
\hat{\mathbb{E}}[\sup_{0\leq s\leq t}|X_s|^p]
\leq&\hat{\mathbb{E}}[\sup_{0\leq s\leq t}|x+\int^s_0 f_u(X_u,K_u)du+\int^s_0 h_u(X_u,K_u)d\langle B \rangle_u +\int^s_0 g_u(X_u,K_u) dB_u+K_s|^p]\\[6pt]
\leq& C(|x|^p+\hat{\mathbb{E}}[\sup_{0\leq s\leq t}|\int^s_0 f_u(X_u,K_u)du|^p]
+\hat{\mathbb{E}}[\sup_{0\leq s\leq t}|\int^s_0 h_u(X_u,K_u)d\langle B \rangle_u|^p]\\[6pt]
+&\hat{\mathbb{E}}[\sup_{0\leq s\leq t}|\int^s_0 g_u(X_u,K_u)dB_u|^p]+\hat{\mathbb{E}}[|K_t|^p]).\end{array}
\end{equation*}
Similarly, from the representation of $K$ (\ref{represen}), we have
\begin{equation*}
\begin{array}{r@{\ }l}
\hat{\mathbb{E}}[K^p_t]\leq&\hat{\mathbb{E}}[\sup_{0\leq s\leq t}(( x+\int^s_0 f_u(X_u,K_u)du+\int^s_0 h_u(X_u,K_u)d\langle B \rangle_u +\int^s_0 g_u(X_u,K_u) dB_u-S_s)^-)^p]\\[6pt]
\leq&\hat{\mathbb{E}}[\sup_{0\leq s\leq t}(( x+\int^s_0 f_u(X_u,K_u)du+\int^s_0 h_u(X_u,K_u)d\langle B \rangle_u +\int^s_0 g_u(X_u,K_u) dB_u-S^+_s)^-)^p]\\[6pt]
\leq&\hat{\mathbb{E}}[\sup_{0\leq s\leq t}|x+\int^s_0 f_u(X_u,K_u)du+\int^s_0 h_u(X_u,K_u)d\langle B \rangle_u +\int^s_0 g_u(X_u,K_u) dB_u-S^+_s|^p]\\[6pt]
\leq& C(|x|^p+\hat{\mathbb{E}}[\sup_{0\leq s\leq t}|\int^s_0 f_u(X_u,K_u)du|^p]
+\hat{\mathbb{E}}[\sup_{0\leq s\leq t}|\int^s_0 h_u(X_u,K_u)d\langle B \rangle_u|^p]\\[6pt]
+&\hat{\mathbb{E}}[\sup_{0\leq s\leq t}|\int^s_0 g_u(X_u,K_u)dB_u|^p]+\hat{\mathbb{E}}[\sup_{0\leq s\leq t}|S^+_s|^p]).
\end{array}
\end{equation*}
Combining the above two inequalities and applying BDG type inequalities, we get
\begin{equation*}
\begin{array}{r@{\ }l}
\hat{\mathbb{E}}[\sup_{0\leq s\leq t}|X_s|^p]+\hat{\mathbb{E}}[K^p_t]\leq C(|x|^p+&\int^t_0 (\hat{\mathbb{E}}[|f_s(X_s,K_s)|^p]\\[6pt]
+&\hat{\mathbb{E}}[|h_s(X_s,K_s)|^p]+\hat{\mathbb{E}}[|g_s(X_s,K_s)|^p])ds+\hat{\mathbb{E}}[\sup_{0\leq s\leq t}|S^+_s|^p]).
\end{array}
\end{equation*}
By condition (A2), we deduce
\begin{equation*}
\begin{array}{r@{\ }l}
\hat{\mathbb{E}}[\sup_{0\leq s\leq t}|X_s|^p]+\hat{\mathbb{E}}[K^p_t]
\leq& C(|x|^p
+\int^t_0\hat{\mathbb{E}}[|\beta_1(s)|^p+\beta^p_2(|X_s|^p+|K_s|^p)]ds+\hat{\mathbb{E}}[\sup_{0\leq s\leq t}|S^+_s|^p]\\[6pt]
\leq&C(|x|^p
+\int^T_0 \hat{\mathbb{E}}[|\beta_1(s)|^p]dt+\hat{\mathbb{E}}[\sup_{0\leq t\leq T}|S^+_t|^p]\\[6pt]
+&\beta^p_2\int^t_0\hat{\mathbb{E}}[\sup_{0\leq u\leq s}|X_u|^p]+\hat{\mathbb{E}}[|K_s|^p]ds).
\end{array}
\end{equation*}
Applying Gronwall's lemma to $\hat{\mathbb{E}}[\sup_{0\leq s\leq t}|X_s|^p]+\hat{\mathbb{E}}[|K_t|^p]$, the result follows.\hfill{}$\square$\\[6pt]
\begin{lemma}\label{345}
Assume that $(x^i, f^i, h^i, g^i, S^i)$ satisfy (A1)-(A4), and let $(X^i, K^i)$ be the solution to the RGSDE corresponding to $(x^i, f^i, h^i, g^i, S^i)$, $i=1, 2$. Define
\begin{align*}
\Delta x = x^1-x^2,\ \Delta f = f^1-f^2,\ &\Delta h = h^1-h^2,\ \Delta g = g^1-g^2;\\
\Delta S = S^1-S^2,\ \Delta X = X^1&-X^2,\ \Delta K = K^1-K^2.
\end{align*}
Then either $\hat{\mathbb{E}}[\sup_{0\leq t\leq T}|\Delta X_t|^p]+\hat{\mathbb{E}}[\sup_{0\leq s\leq T}|\Delta K_t|^p]=0$ or there exists a constant $C>0$ such that
\begin{equation*}
\begin{array}{r@{\ }l}
\hat{\mathbb{E}}[\sup_{0\leq t\leq T}|\Delta X_t|^p]+\hat{\mathbb{E}}[\sup_{0\leq s\leq T}|\Delta K_t|^p]\leq \upsilon^{-1}(\upsilon(C(|\Delta x|^p
+&\int^T_0(\hat{\mathbb{E}}[|\Delta f_t(X_t^1,K_t^1)|^p]
+\hat{\mathbb{E}}[|\Delta h_t(X_t^1,K_t^1)|^p]\\
+&\hat{\mathbb{E}}[|\Delta g_t(X_t^1,K_t^1)|^p])dt+\hat{\mathbb{E}}[\sup_{0\leq t\leq T}|\Delta S_t|^p]))\\
+& C\int^T_0{\beta(t)dt}).
\end{array}
\end{equation*}
\end{lemma}
\noindent\textbf{Proof :} We set
$$
\begin{array}{r@{\ }l}
(M^X)^i_t=x^i+\int^t_0 f^i_s(X^i_s,K^i_s)ds+\int^t_0 h^i_s(X^i_s,K^i_s)d\langle B \rangle_s +\int^t_0 g^i_s(X^i_s,K^i_s) dB_s,\ 0\leq t\leq T,\ i=1, 2,\\
\end{array}
$$
and $\Delta M^X=(M^X)^1-(M^X)^2$. By condition (A3), we have
$$
\begin{array}{r@{\ }l}
\hat{\mathbb{E}}[\sup_{0\leq s\leq t}|(\Delta M^X)_s|^p]\leq& \hat{\mathbb{E}}[\sup_{0\leq s\leq t}|\Delta x
+\int^s_0 (f^1_u(X^1_u,K^1_u)-f^2_u(X^2_u,K^2_u))du\\[6pt]
+&\int^s_0 (h^1_u(X^1_u,K^1_u)-h^2_u(X^2_u,K^2_u))d\langle B\rangle_u
+\int^s_0 (g^1_u(X^1_u,K^1_u)-g^2_u(X^2_u,K^2_u))dB_u|^p]\\[6pt]
\leq&\hat{\mathbb{E}}[\sup_{0\leq s\leq t}|\Delta x
+\int^s_0 \Delta f_u(X_u^1,K_u^1) du
+\int^s_0 (f^2_u(X^1_u,K^1_u)-f^2_u(X^2_u,K^2_u))du\\[6pt]
+&\int^s_0 \Delta h_u(X_u^1,K_u^1)d\langle B\rangle_u
+\int^s_0 (h^2_u(X^1_u,K^1_u)-h^2_u(X^2_u,K^2_u))d\langle B\rangle_u\\[6pt]
+&\int^s_0 \Delta g_u(X_u^1,K_u^1)dB_u
+\int^s_0 (g^2_u(X^1_u,K^1_u)-g^2_u(X^2_u,K^2_u))dB_u|^p]\\[6pt]
\leq& C(|\Delta x|^p
+\int^t_0(\hat{\mathbb{E}}[|\Delta f_s(X_s^1,K_s^1)|^p]
+\hat{\mathbb{E}}[|\Delta h_s(X_s^1,K_s^1)|^p]\\[6pt]+&\hat{\mathbb{E}}[|\Delta g_s(X_s^1,K_s^1)|^p])ds
+\int^t_0\beta(s)\rho(\hat{\mathbb{E}}[|\Delta X_s|^p]+\hat{\mathbb{E}}[|\Delta K_s|^p])ds),
\end{array}
$$
and
\begin{equation*}
\begin{array}{r@{\ }l}
\hat{\mathbb{E}}[\sup_{0\leq s\leq t}|\Delta K_s|^p]&=\hat{\mathbb{E}}[\sup_{0\leq s\leq t}|\sup_{0\leq u\leq s}((M^X)^1_u-S^1_u)^--\sup_{0\leq u\leq s}((M^X)^2_u-S^2_u)^-|^p]\\[6pt]
&\leq\hat{\mathbb{E}}[\sup_{0\leq s\leq t}|\sup_{0\leq u\leq s}|{((M^X)^1_u-S^1_u)^-}-((M^X)^2_u-S^2_u)^-||^p]\\[6pt]
&=\hat{\mathbb{E}}[\sup_{0\leq s\leq t}|{((M^X)^1_s-S^1_s)^-}-((M^X)^2_s-S^2_s)^-|^p]\\[6pt]
&\leq\hat{\mathbb{E}}[\sup_{0\leq s\leq t}|((M^X)^1_s-S^1_s)-((M^X)^2_s-S^2_s)|^p]\\[6pt]
&\leq C(\hat{\mathbb{E}}[\sup_{0\leq s\leq t}|\Delta (M^X)_s|^p]+\hat{\mathbb{E}}[\sup_{0\leq s\leq t}|\Delta S_s|^p]).
\end{array}
\end{equation*}
Then, we have
$$
\begin{array}{r@{\ }l}
\hat{\mathbb{E}}[\sup_{0\leq s \leq t}|\Delta X_s|^p]+\hat{\mathbb{E}}[\sup_{0\leq s \leq t}|\Delta K_s|^p]
\leq& \hat{\mathbb{E}}[\sup_{0\leq s\leq t}|(\Delta M^X)_s+\Delta K_s|^p]+\hat{\mathbb{E}}[\sup_{0\leq s \leq t}|\Delta K_s|^p]\\[6pt]
\leq& C(\hat{\mathbb{E}}[\sup_{0\leq s\leq t}|(\Delta M^X)_s|^p]+\hat{\mathbb{E}}[\sup_{0\leq s\leq t}|\Delta K_s|^p])\\[6pt]
\leq& C(|\Delta x|^p
+\int^t_0(\bar{\mathbb{E}}[|\Delta f_s(X_s^1,K_s^1)|^p]
+\hat{\mathbb{E}}[|\Delta h_s(X_s^1,K_s^1)|^p]\\[6pt]+&\hat{\mathbb{E}}[|\Delta g_s(X_s^1,K^s_1)|^p])ds
+\hat{\mathbb{E}}[\sup_{0\leq s\leq t}|\Delta S_s|^p]\\+&\int^t_0\beta(s)\rho(\bar{\mathbb{E}}[|\Delta X_s|^p+|\Delta K_s|^p])ds).\\[6pt]
\leq& C(|\Delta x|^p
+\int^T_0(\bar{\mathbb{E}}[|\Delta f_t(X_t^1,K_t^1)|^p]
+\hat{\mathbb{E}}[|\Delta h_t(X_t^1,K_t^1)|^p]\\+&\hat{\mathbb{E}}[|\Delta g_t(X_t^1,K_t^1)|^p])dt+\hat{\mathbb{E}}[\sup_{0\leq t\leq T}|\Delta S_t|^p]\\
+&\int^t_0\beta(s)\rho(\hat{\mathbb{E}}[\sup_{0\leq u\leq s}|\Delta X_u|^p]+\bar{\mathbb{E}}[\sup_{0\leq u\leq s}|\Delta K_u|^p])ds).\\[6pt]
\end{array}
$$
Lemma \ref{le 2.2.1} gives the desired result.\hfill{}$\square$\\[6pt]
Now we will give the proof of Theorem \ref{Exi}.\\
\noindent\textbf{Proof of Theorem \ref{Exi}:}
By taking $x^1=x^2$, $f^1=f^2$, $h^1=h^2$, $g^1=g^2$ and $S^1=S^2$ in Lemma \ref{345}, we obtain immediately the uniqueness result. Now we will establish the existence result for RGSDE with nonlinear resistance (\ref{1}) by using a Picard iteration. Letting $X^0=x$ and $K^0=0$, for each $n\in\mathbb{N}_+$, $X^{n+1}$ and $K^{n+1}$ are given by recurrence:
\begin{align*}
&X^{n+1}_t=x+\int^t_0 f_s(X^n_s,K^n_s)ds+\int^t_0 h_s(X^n_s,K^n_s)d\langle B \rangle_s +\int^t_0 g_s(X^n_s,K^n_s) dB_s + K^{n+1}_t,\ 0\leq t\leq T,
\end{align*}
satisfying\\[6pt]
$
\begin{array}{l}
(a)\ X^{n+1}\in M^p_G([0, T]),\ X^{n+1}_t\geq S_t,\ q.s.;\\[6pt]
(b)\ K^{n+1}\in M_I([0, T])\cap M^p_G([0, T]),\ K^{n+1}_0=0,\ q.s.;\\[6pt]
(c)\ \int^T_0(X_t^{n+1}-S_t)dK^{n+1}_t=0,\ q.s..
\end{array}
$\\[6pt]
Actually $X^{n+1}_t=\tilde{X}^{n+1}_t$ and $K^{n+1}_t=\tilde{K}^{n+1}_t$, where $(\tilde{X}^{n+1},\tilde{K}^{n+1})$ is given by Lemma \ref{BXT} with $Y_t=x+\int^t_0 f_s(X^n_s,K^n_s)ds+\int^t_0 h_s(X^n_s,K^n_s)d\langle B \rangle_s +\int^t_0 g_s(X^n_s,K^n_s) dB_s-S_t$. Thus $(X^{n+1}, K^{n+1})$ is well defined in $M^p_G([0,T])\times (M_I([0, T])\cap M^p_G([0,T]))$.\\[6pt]
Similarly to Lemma \ref{abcde}, we get a priori estimate uniform in $n$ for $\{\hat{\mathbb{E}}[\sup_{0\leq t\leq T}|X^n_t|^p]+\hat{\mathbb{E}}[|K^{n}_T|^p]\}_{n\in \mathbb{N}}$. Indeed, we have
\begin{align*}
\hat{\mathbb{E}}[\sup_{0\leq s\leq t}|X^{n+1}_t|^p]+\hat{\mathbb{E}}[|K^{n+1}_t|^p]
\leq C\bigg(|x|^p
+&\int^t_0 \hat{\mathbb{E}}[|\beta_1(s)|^p]ds+\hat{\mathbb{E}}[\sup_{0\leq s\leq t}|S^+_s|^p]\\
+&\int^t_0 \beta^p_2(\hat{\mathbb{E}}[\sup_{0\leq u\leq s}|X^n_u|^p]+\hat{\mathbb{E}}[|K^n_s|^p])ds\bigg).
\end{align*}
Set
\[
a(t):=Ce^{C\beta^p_2t}\bigg(|x|^p+\int^T_0 \hat{\mathbb{E}}[|\beta_1(s)|^p]ds+\hat{\mathbb{E}}[\sup_{0\leq s\leq T}|S^+_s|^p]\bigg),
\]
then $a(\cdot)$ is the solution to the following ordinary differential equation:
\[
a(t)=C\bigg(|x|^p+\int^T_0 \hat{\mathbb{E}}[|\beta_1(s)|^p]ds+\hat{\mathbb{E}}[\sup_{0\leq s\leq T}|S^+_s|^p]+\beta^p_2\int^t_0a(s)ds\bigg).
\]
It is easy to check that for all $n\in\mathbb{N}$, $\hat{\mathbb{E}}[\sup_{0\leq s\leq t}|X^{n}_t|^p]+\hat{\mathbb{E}}[|K^{n}_t|^p]\leq a(t)$.\\[6pt]
Secondly, for $n$ and $m\in \mathbb{N}$, we define
$$u^{n+1,m}_t:=\hat{\mathbb{E}}[\sup_{0\leq s \leq t}|X^{n+m+1}_s-X^{n+1}_s|^p]+\hat{\mathbb{E}}[\sup_{0\leq s \leq t}|K^{n+m+1}_s-K^{n+1}_s|^p],\ 0\leq t\leq T.$$
Then it follows from the same techniques of the proof of Theorem 3.7 in Lin and Bai \cite{BL} that $\{X^n\}_{n\in\mathbb{N}}$ and $\{K^n\}_{n\in\mathbb{N}}$ are two Cauchy sequences in $M^p_G([0, T])$. We denote the limit by $X$ and $K$. Following the procedures of the proof of Lemma \ref{345} and noting that $\rho$ is continuous and $\rho(0_+)=0$, $K$ has the representation (\ref{represen}). \\
Obviously, the pair of processes $(X, K)$ satisfies (A5) - (A7). Thus the pair of processes $(X,K)$, well defined in $M^p_G([0, T])\times (M_I([0, T])\cap M^p_G([0, T]))$, is a solution to (\ref{1}).\hfill{}$\square$\\[6pt]

In order to give the comparison principle, we consider the following scalar valued RGSDE:
\begin{equation*}
X_t=x+\int^t_0 f_s(X_s,K_s)ds+\int^t_0 h_s(X_s,K_s)d\langle B \rangle_s +\int^t_0 g_s(X_s) dB_s +K_t,\ 0\leq t\leq T,
\end{equation*} and we assume that:
\begin{mnlist}\addtocounter{mnelist}{1}
%\item the initial condition $x\in \mathbb{r}$;
\item For some $p>2$, the coefficients $f$, $h$: $\Omega\times[0,T]\times\mathbb{R}\times\mathbb{R}\rightarrow\mathbb{R}$ and $g:\Omega\times[0,T]\times\mathbb{R}\rightarrow\mathbb{R}$ are given functions satisfying for each $x,y\in\mathbb{R}$, $f_\cdot(x,y)$, $h_\cdot(x,y)$, and $g_\cdot(x)\in M^p_G([0,T])$ and $|f_t(x,y)|^p+|h_t(x,y)|^p+|g_t(x)|^p\leq |\beta_1(t)|^p+\beta^{p}_{2}(|x|^p+|y|^p),$ where $\beta_1\in M^p_G([0,T])$ and $\beta_2\in\mathbb{R}_{+}$;
\item The coefficients $f$, $h$ and $g$ satisfying an integral-Lipschitz condition, i.e., for each $t\in[0, T]$ and $x, x', y, y'\in\mathbb{R}$, $|f_t(x,y)-f_t(x',y')|^p+|h_t(x,y)-h_t(x',y')|^p+|g_t(x)-g_t(x')|^p\leq \rho(|x-x'|^p+|y-y'|^p),$
    where $\rho~:(0,+\infty)\rightarrow(0,+\infty)$ is  continuous increasing and concave function that vanishes at $0_+$ and satisfies
 \[
 \int^1_0{\frac{dr}{r+\rho(r)}}=+\infty.
 \]

% where $\beta~:[0,T]\rightarrow\mathbb{R}^{+}$ is integrable, and $\rho~:(0,+\infty)\rightarrow(0,+\infty)$ is  continuous increasing and concave function satisties
% \[
% \rho(0)=0,~~\int^1_0{\frac{dr}{\rho(r)}}=+\infty
% \]
%\item The obstacle is a $G$-It\^o process whose coefficients are all elements in $M^p_G([0,T])$,
%and we shall always assume that $S_0\leq x$, q.s..
\end{mnlist}
 \begin{remark}
 Before we move to the comparison principle, we should mention that restricted to the fact that we need to apply $G$-It\^o's formula and we need to consider $X$ separately, we should impose stronger conditions on the coefficients but which are still weaker than those assumed in Lin \cite{L}. It is easy to check that (A3') implies (A3). Following the comparison result in Lin \cite{L}, at first, we assume that coefficients $f$, $h$ and $g$ and the obstacle process $S$ are bounded, and then we remove it in the second step.
 \end{remark}
We then have the following results.
\begin{theorem}\label{CT}
Suppose that for $i=1,2$, $f^i,~h^i,~g^i$ satisfy the conditions (A1),(A2'),(A3') and (A4), and we assume the following:
\begin{mylist}
\item $x^1\leq x^2$;
\item $f^i$, $h^i$ and $g^1=g^2=g$ are bounded, and $S^i$ are uniformly upper bounded, $i=1, 2$;
\item $f^1_t(x,0)\leq f^2_t(x,0)$ and $h^1_t(x,0)\leq h^2_t(x,0)$, for $x\in\mathbb{R}$, $f^1,h^1$ are decreasing in $y$, and $f^2,h^2$ are increasing in $y$, and $S^1_t\leq S^2_t$, $0\leq t\leq T$, q.s..
\end{mylist}
If $(X^i, K^i)$ is the solution to the RGSDEs with data $(f^i, h^i, g, S^i)$, $i=1, 2$, then,
$$X^1_t\leq X^2_t,\ 0\leq t\leq T,\ q.s..$$
\end{theorem}
\noindent\textbf{Proof :}
Since $f^i$, $h^i$ and $g$ are bounded, and $S^i$ are uniformly upper bounded, $i=1, 2$,
using the BDG type inequalities to (\ref{represen}), we deduce that $K^i_T$ has the moment for arbitrage large order and for $0\leq t\leq T$, $\lim_{s\rightarrow t}\bar{\mathbb{E}}[|K^i_t-K^i_s|^2]=0$, $i=1, 2$. \\[6pt]
\noindent Compared to Lin \cite{L}, we consider $(x^+)^p$ and apply the extended $G$-It\^o's formula to
$((X^1_t-X^2_t)^+)^p$,
\begin{align}\label{cal1}
((X^1_t-X^2_t)^+)^p
=&3\int^t_0|(X^1_s-X^2_s)^+|^{p-1}(f^1_s(X^1_s,K^1_s)-f^2_s(X^2_s,K^2_s))ds\notag\\
+&3\int^t_0|(X^1_s-X^2_s)^+|^{p-1}(h^1_s(X^1_s,K^1_s)-h^2_s(X^2_s,K^2_s))d\langle B\rangle_s\notag\\
+&3\int^t_0|(X^1_s-X^2_s)^+|^{p-1}(g_s(X^1_s)-g_s(X^2_s))dB_s\\
+&3\int^t_0|(X^1_s-X^2_s)^+|^{p-1}d(K^1_s-K^2_s)\notag\\
+&3\int^t_0|(X^1_s-X^2_s)^+|^{p-2}|g_s(X^1)-g_s(X^2_s)|^2 d\langle B\rangle_s\notag
\end{align}
Since on $\{X^1_t>X^2_t\}$, $X^1_t>X^2_t\geq S^2_t\geq S^1_t$, we have
\begin{align}\label{cal2}
\int^t_0|(X^1_s-X^2_s)^+|^{p-1}d(K^1_s-K^2_s)&=\int^t_0|(X^1_s-X^2_s)^+|^{p-1}dK^1_s-
\int^t_0|(X^1_s-X^2_s)^+|^{p-1}dK^2_s\notag\\
&\leq \int^t_0|(X^1_s-S^1_s)^+|^{p-1}dK^1_s-
\int^t_0|(X^1_s-X^2_s)^+|^2dK^{p-1}_s\\
&\leq -
\int^t_0|(X^1_s-X^2_s)^+|^{p-1}dK^2_s\leq 0,\ q.s..\notag
\end{align}
Noting the monotonicity of $f^1$, $h^1$, $f^2$, and $h^2$ and (\ref{cal2}), we have
\begin{align}\label{ca}
((X^1_t-X^2_t)^+)^p
\leq&3\int^t_0|(X^1_s-X^2_s)^+|^{p-1}(f^1_s(X^1_s,0)-f^2_s(X^2_s,0))ds\notag\\
+&3\int^t_0|(X^1_s-X^2_s)^+|^{p-1}(h^1_s(X^1_s,0)-h^2_s(X^2_s,0))d\langle B\rangle_s\notag\\
+&3\int^t_0|(X^1_s-X^2_s)^+|^{p-1}(g_s(X^1_s)-g_s(X^2_s))dB_s\\
+&3\int^t_0|(X^1_s-X^2_s)^+|^{p-2}|g_s(X^1)-g_s(X^2_s)|^2 d\langle B\rangle_s\notag
\end{align}
Then, by condition (A3'), Young's inequality, taking $G$-expectation on both sides of (\ref{ca}) and Jensen's inequality, we conclude
\begin{align*}\label{complex}
\hat{\mathbb{E}}[((X^1_t-X^2_t)^+)^p]&\leq C\hat{\mathbb{E}}\bigg[\int^t_0((X^1_s-X^2_s)^+)^p+\rho(((X^1_s-X^2_s)^+)^p)ds\bigg]\\
&\leq C\int^t_0\hat{\mathbb{E}}[((X^1_s-X^2_s)^+)^p]+\rho(\hat{\mathbb{E}}[((X^1_s-X^2_s)^+)^p])ds.
\end{align*}
Using Bihari¡¯s inequality, it follows that $\hat{\mathbb{E}}[((X^1_t-X^2_t)^+)^p]=0$, which implies the result.\hfill{}$\square$\\[6pt]
\begin{theorem}\label{C}
Suppose that for $i=1,2$ $f^i,~h^i,~h^i$ satisfy the conditions (A1),(A2'),(A3') and (A4), and we assume in addition the following:
\begin{mylist}
\item $x^1\leq x^2$ and $g^1=g^2=g$;
\item $f^1_t(x,0)\leq f^2_t(x,0)$ and $h^1_t(x,0)\leq h^2_t(x,0)$, for $x\in\mathbb{R}$, $f^1,h^1$ are decreasing in $y$, and $f^2,h^2$ are increasing in $y$, and $S^1_t\leq S^2_t$, $0\leq t\leq T$, q.s..
\end{mylist}
If $(X^1,K^1)$ and $(X^2, K^2)$ are the solutions to the RGSDEs above resepectively, then,
$$X^1_t\leq X^2_t,\ 0\leq t\leq T,\ q.s..$$
\end{theorem}
\noindent\textbf{Proof :} Firstly, we define the truncation functions for the coefficients and the obstacle process: for $N>0$, $\xi^N_t(x)=(-N\vee \xi_t(x))\wedge N$, $x\in\mathbb{R}$, where $\xi$ denote $f^i$, $h^i$ and $g$, $i=1, 2$, and $S^N_t=S_t\wedge N$, $0\leq t\leq T$.
It is easy to verify that the truncated coefficients and obstacle processes satisfy (A2') and (A3').
Moreover, the truncation functions keep the order of the coefficients and obstacle processes, that is,
$$(f^1)^N_t(x,0)\leq (f^2)^N_t(x,0),\ (h^1)^N_t(x,0)\leq (h^2)^N_t(x,0), \text{ and } (S^1)^N_t \leq (S^2)^N_t,\ 0\leq t\leq T,\ q.s..$$
$$(f^1)^N_t,(h^1)^N_t \text{ are decreasing in }y,\text{ and }\ (f^2)^N_t,(h^2)^N_t \text{ are increasing in } y.$$
Consider the following RGSDEs on $[0, T]$, for $i=1, 2$,
\begin{align*}
(X^i)^N_t=x+\int^t_0 (f^i)^N_s((X^i)^N_s,(K^i)^N_s)ds+\int^t_0 (h^i)^N_s((X^i)^N_s,(K^i)^N_s)d\langle B \rangle_s +\int^t_0 g^N_s((X^i)^N_s) dB_s + (K^i)^N_t;
\end{align*}
satisfies\\[6pt]
$
\begin{array}{l}
(a)\ (X^i)^N_t\in M^p_G([0, T]),\ (X^i)^{N}_t\geq (S^i)^N_t,\ q.s.;\\[6pt]
(b)\ (K^i)^N\in M_I([0, T])\cap M^p_G([0, T]),\ (K^i)^N_0=0,\ q.s.;\\[6pt]
(c)\ \int^T_0((X^i)^N_t-(S^i)^N_t)d(K^i)^{N}_t=0,\ q.s..
\end{array}
$\\[6pt]
By Theorem \ref{CT},
we have
\begin{equation}\label{last}
(X^1)^N_t\leq (X^2)^N_t,\ q.s..
\end{equation}
In view of the proof of Lemma \ref{345}, we have
%$$
%\begin{array}{r@{\ }l}
%\bar{\mathbb{E}}[\sup_{0\leq s\leq t}|(X^i)^N_s-X^i_s|^p]
%\leq&C(\int^T_0(
%\bar{\mathbb{E}}[|(f^i)^N_t(X^i_t)-f^i_t(X^i_t)|^p]+\bar{\mathbb{E}}[|(h^i)^N_t(X^i_t)-h^i_t(X^i_t)|^p]\\[6pt]
%+&\bar{\mathbb{E}}[|g^N_t(X^i_t)-g_t(X^i_t)|^p]
%)dt+\bar{\mathbb{E}}[\sup_{0\leq t\leq T}|(S^i)^N_t-S^i_t|^p]\\[6pt]
%+&\int^t_0\bar{\mathbb{E}}[\sup_{0\leq u\leq s}|(X^i)^N_u-X^i_u|^p]ds).
%\end{array}
%$$
%Applying again Gronwall's lemma, we obtain
%\begin{align*}
%\bar{\mathbb{E}}[\sup_{0\leq t\leq T}|(X^i)^N_t-X^i_t|^p]
%\leq& C\bigg(\int^T_0(
%\bar{\mathbb{E}}[|(f^i)^N(t, X^i_t,K^i_t)-f^i(t, X^i_t,K^i_t)|^p]+\bar{\mathbb{E}}[|(h^i)^N(t, X^i_t,K^i_t)-h^i(t, X^i_t,K^i_t)|^p]\\
%+&\bar{\mathbb{E}}[|g^N(t, X^i_t)-g(t, X^i_t)|^p]
%)dt+\bar{\mathbb{E}}[\sup_{0\leq t\leq T}|(S^i)^N_t-S^i_t|^p]\bigg).
%\end{align*}
\begin{equation*}
\begin{array}{r@{\ }l}
\hat{\mathbb{E}}[\sup_{0\leq t\leq T}|(X^i)^N_t-X^i_t|^p]+\hat{\mathbb{E}}[\sup_{0\leq t\leq T}|(K^i)^N_t-K^i_t|^p]&\leq C(\int^T_0(\hat{\mathbb{E}}[|(f^i)^N(t, X^i_t,K^i_t)-f^i(t, X^i_t,K^i_t)|^p]\\
+&\hat{\mathbb{E}}[|(h^i)^N(t, X^i_t,K^i_t)-h^i(t, X^i_t,K^i_t)|^p]\\
+&\hat{\mathbb{E}}[|(g)^N(t, X^i_t)-g(t, X^i_t)|^p])dt\\
+&\hat{\mathbb{E}}[\sup_{0\leq t\leq T}|(S^i)^N_t-S^i_t|^p])).
\end{array}
\end{equation*}
For any $t\in[0, T]$, by condition (A2') we calculate
\begin{align*}
\hat{\mathbb{E}}[|(f^i)^N_t(X^i_t,K^i_t)-f^i_t(X^i_t,K^i_t)|^p]\leq& \hat{\mathbb{E}}[|f^i_t(X^i_t,K^i_t)|^pI_{\{|f^i_t(X^i_t,K^i_t)|>N\}}]\\
\leq& \hat{\mathbb{E}}[(|\beta_1(t)|^p+\beta^{p}_2(|X^i_t|^p+K|^i_t|^p))I_{\{(|\beta_1(t)|^p+\beta^{p}_2(|X^i_t|^p+K|^i_t|^p))>N^p\}}]\\
\leq& C(\hat{\mathbb{E}}[|\beta_1(t)|^pI_{\{|\beta_1(t))|^p>\frac{N^p}{3}\}}]+\hat{\mathbb{E}}[|X^i_t|^pI_{\{|X^i_t|^p>\frac{N^p}{3\beta^{p}_2}\}}]+\bar{\mathbb{E}}[|K^i_t|^pI_{\{|K^i_t|^p>\frac{N^p}{3\beta^{p}_2}\}}]).
\end{align*}
Since $\beta_1(\cdot)$, $X^i$ and $K^i\in M^p_G([0, T])$, we obtain we that $\beta_1(t)$, $X^i_t$ and $K^i_t\in L^p_G(\Omega_t)$ for almost every $t\in[0, T]$. Therefore, letting $N\rightarrow +\infty$, we have
$$
\hat{\mathbb{E}}[|(f^i)^N_t(X^i_t,K^i_t)-f^i_t(X^i_t,K^i_t)|^p]\rightarrow 0.
$$
Similarly, we can obtain that
$$
\hat{\mathbb{E}}[|(h^i)^N_t(X^i_t,K^i_t)-h^i_t(X^i_t,K^i_t)|^p]\rightarrow 0
;
$$
and
$$
\hat{\mathbb{E}}[|(g^i)^N_t(X^i_t)-g^i_t(X^i_t)|^p]\rightarrow 0.
$$
%\noindent Using dominated convergence theorem to the integrals on $[0, T]$, it follows that
%\begin{align}\label{cal3}
%\lim_{N\rightarrow +\infty}
%\int^T_0(\bar{\mathbb{E}}[|(f^i)^N(t, X^i_t,K^i_t)-f^i(t, X^i_t,K^i_t)|^p]+\bar{\mathbb{E}}[|(h^i)^N(t, X^i_t,K^i_t)-&h^i(t, X^i_t,K^i_t)|^p]\\
%+&\bar{\mathbb{E}}[|g^N(t, X^i_t)-g(t, X^i_t)|^p]
%)dt=0.\notag
%\end{align}
%On the other hand,
%$$
%\bar{\mathbb{E}}[\sup_{0\leq t\leq T}|(S^i)^N_t-S^i_t|^p]\leq \bar{\mathbb{E}}[\sup_{0\leq t\leq T}(|S^i_t|^pI_{\{|S^i_t|>N\}})]\leq \bar{\mathbb{E}}[\sup_{0\leq t\leq T}|S^i_t|^pI_{\{\sup_{0\leq t\leq T}|S^i_t|>N\}}].
%$$
%By the proof of Theorem \ref{BXT}, we know that $\sup_{0\leq t\leq T}S^i_t$ is an element in $L^p_G(\Omega_T)$. So we have
%\begin{align}\label{cal4}
%\bar{\mathbb{E}}[\sup_{0\leq t\leq T}|(S^i)^N_t-S^i_t|^p]\leq \bar{\mathbb{E}}[\sup_{0\leq t\leq T}|S^i_t|^pI_{\{\sup_{0\leq t\leq T}|S^i_t|>N\}}]\rightarrow 0,\ as\ N\rightarrow +\infty.
%\end{align}
Thus, it follows from the procedures of the proof of Theorem 5.9 in Lin \cite{L} that
\begin{equation}\label{last2}
\hat{\mathbb{E}}[\sup_{0\leq t\leq T}|(X^i)^N_t-X^i_t|^p]\rightarrow 0,\ as\ N\rightarrow +\infty.
\end{equation}
Then, (\ref{last}) and (\ref{last2}) yield the desired result .\hfill{}$\square$\\[6pt]
As a consequent result of Theorem \ref{C}, we have the following corollary.
\begin{corollary}
Given two RGSDEs satisfying the conditions (A1),(A2'),(A3') and (A4), and suppose that $x^1\leq x^2$ and $g^1=g^2=g$, if one of the following holds:
\begin{mylist}
\item $f^1$ and $h^1$ are independent of y, and $f^1_t(x)\leq f^2_t(x,0)$ and $h^1_t(x)\leq h^2_t(x,0)$, for $x\in\mathbb{R}$, $f^2,h^2$ are increasing in $y$, and $S^1_t\leq S^2_t$, $0\leq t\leq T$, q.s.;
\item $f^2$ and $h^2$ are independent of y, and $f^1_t(x,0)\leq f^2_t(x)$ and $h^1_t(x,0)\leq h^2_t(x)$, for $x\in\mathbb{R}$, $f^1,h^1$ are decreasing in $y$, and $S^1_t\leq S^2_t$, $0\leq t\leq T$, q.s..
\end{mylist}
Let $(X^1,K^1)$ and $(X^2, K^2)$ are two pairs of solutions to the RGSDEs above respectively, then,
$$X^1_t\leq X^2_t,\ 0\leq t\leq T,\ q.s..$$
\end{corollary}

%we get
%\begin{equation}\label{est4}
%\begin{array}{r@{\ }l}
%\bar{\mathbb{E}}[\sup_{0\leq s\leq t}|X_t|^p]+\bar{\mathbb{E}}[|K_t|^p]\leq \upsilon^{-1}(\upsilon(C(|x|^p
%+&\int^T_0(\bar{\mathbb{E}}[|f_t(0)|^p]
%+\bar{\mathbb{E}}[|h_t(0)|^p]\\[6pt]
%+&\bar{\mathbb{E}}[|g_t(0)|^p])dt
%+\bar{\mathbb{E}}[\sup_{0\leq t\leq T}|S^+_t|^p])+C\int^t_0{\beta(s)ds}),\ 0\leq t\leq T.
%\end{array}
%\end{equation}
%Putting (\ref{est4}) into (\ref{est3}),

\noindent\textbf{Acknowledgement}\ %The author expresses special thanks to Prof. Ying HU, who provided both the initial inspiration for the work and useful suggestions, and also to the anonymous reviewers, who have given me a lot of important constructive advice for the revision.

\end{document}